\documentclass[a4paper,11pt]{article}

\usepackage{CAR}
\usepackage[draft=true]{minted}
\usepackage{booktabs}
\usepackage{cleveref}
\usepackage{tcolorbox}
\BeforeBeginEnvironment{minted}{\begin{tcolorbox}}
\AfterEndEnvironment{minted}{\end{tcolorbox}}
\usepackage{cite}

\begin{document}

\setcounter{footnote}{0}
\setcounter{figure}{0}


\Abschnitt
{Neues über Systeme}
{Neues über Systeme}
{rubrik}

\vspace{3mm}


\Aufsatz
{Toric Geometry in OSCAR}
{Toric Geometry in OSCAR}
{Martin Bies, Lars Kastner}
{M. Bies, L. Kastner}
{Martin Bies (RPTU Kaiserslautern-Landau)\\ Lars Kastner (TU Berlin)}
{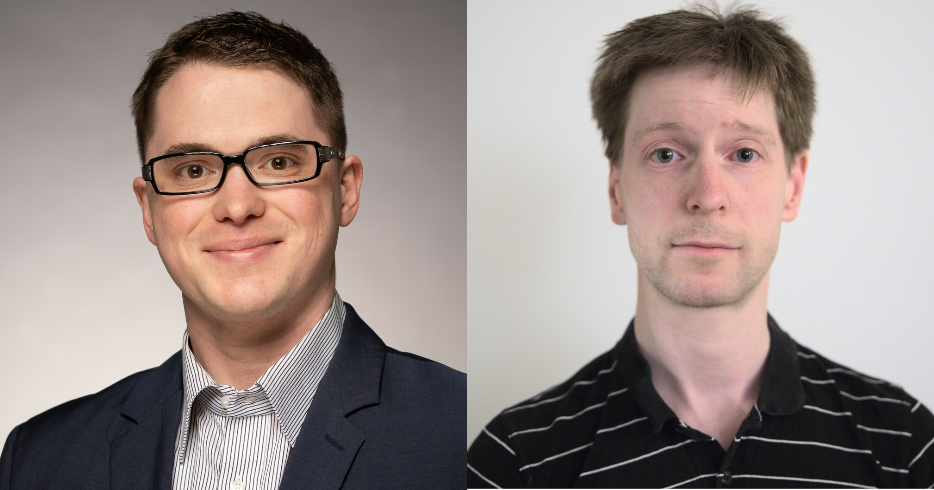}
{bies@mathematik.uni-kl.de,kastner@math.tu-berlin.de}

\begin{otherlanguage}{english}

\vspace{3mm}
\begin{multicols}{2}
\noindent


\Ueberschrift{Why toric geometry in OSCAR?}{sec1}

\Ueberschriftu{Toric geometry -- an arena for mathematical theories}
Among the fields of algebraic geometry, the field of toric geometry is particularly well understood and algorithmic. Among others, the cohomology ring, the Chow ring, topological intersection numbers as well as cohomologies of coherent sheaves can be obtained with computer algorithms \cite{CLS11}. Therefore, toric varieties provide a useful platform for testing mathematical theories.

To put it briefly, toric varieties are characterized by having an algebraic torus $\left( \mathbb{C}^\ast \right)^r$ as a dense and open subset. This is why they are called \emph{toric}. While the realm of toric varieties is more constrained compared to that of general schemes/varieties, the toric universe still provides a significant degree of versatility. As an example, many Calabi-Yau manifolds can be constructed as complete intersections in toric varieties \cite{KS98, KS00}. This includes many K3 surfaces \cite{KS98}. More recently, starting from such K3 surfaces, researchers have discovered in the framework of F-theory -- a non-perturbative regime of string theory -- the largest currently-known class of globally consistent Standard Model solutions without chiral exotics and gauge coupling unification \cite{CHLLT19}.

For all these reasons, there is a high demand for computer implementations of toric geometry. Some examples of computer algebra systems that support toric geometry are \cite{Dev23, GS23}.

\Ueberschriftu{OSCAR -- a melting pot}
We present a computer implementation for toric varieties in the computer algebra system \texttt{OSCAR} \cite{**key**3, **key*}. The funding for \texttt{OSCAR} is provided by the SFB-TRR 195 \emph{Symbolic Tools in Mathematics and their Application} of the German Research Foundation (DFG). The main architectural feature of \texttt{OSCAR} is that its four fundamental tools \texttt{Antic} (Hecke, Nemo), \texttt{GAP}, \texttt{Polymake} and \texttt{Singular} are \emph{integral components}, rather than external software that can be used. For more information, the interested reader can consult the article ``OSCAR:  Open Source Computer Algebra Research system'' by \emph{Prof. Dr. Max Horn} (to appear in the \emph{ComputerAlgebraRundbrief}) or the \texttt{OSCAR} homepage:
\begin{center}
\href{https://www.oscar-system.org}{https://www.oscar-system.org}
\end{center}

By leveraging \texttt{Polymake}, we can carry out polyhedral geometry operations, such as handling cones and fans, and utilize cutting-edge algorithms for triangulations \cite{JJK18}. This provides a reliable backbone for toric geometry in \texttt{OSCAR}. The Cox ring as well as the Chow ring of toric varieties are polynomial rings. Closed subvarieties of toric varieties correspond to homogeneous polynomial in the Cox ring \cite{CLS11}. This functionality is provided by the software \texttt{Singular}. Additionally, tools from group and number theory are essential in toric geometry. Such tasks are executed with \texttt{Antic} (Hecke, Nemo) and \texttt{GAP}. To sum up, toric geometry benefits greatly from the combination of \texttt{Antic} (Hecke, Nemo), \texttt{GAP}, \texttt{Polymake} and \texttt{Singular}.

\Ueberschriftu{\texttt{Julia} -- a modern programming language}
\texttt{Julia} \cite{BEKS17} is a high-performance programming language designed for numerical and scientific computing. The growing ecosystem of \texttt{Julia} packages ensures its continued viability for scientific computing and data analysis. \texttt{OSCAR} is written in \texttt{Julia}. This implies that the performance of \texttt{OSCAR} should be comparable or even better than many other implementations.

\Ueberschrift{Overview}{sec2}

Our goal with \texttt{OSCAR} is to create a computer algebra system that is both user-friendly and convenient. To assist users with toric geometry, we offer a tutorial:\footnote{Interested readers may also explore the actual \texttt{OSCAR} code on \href{https://github.com/oscar-system/Oscar.jl}{GitHub}.}
\begin{center}
\href{https://www.oscar-system.org/tutorials/}{https://www.oscar-system.org/tutorials/}. 
\end{center}
The toric implementation in \texttt{OSCAR} are conceptually based on \cite{CLS11}. This is a fundamental guiding principle within \texttt{OSCAR}: Implementations are conceptually grounded in a few carefully selected publications.

The relationship between toric geometry and polyhedral geometry is crucial for any toric geometry implementation. We illustrate this connection for affine toric varieties.

In \texttt{OSCAR}, the toric implementations focus on the lattice $N = \mathbb{Z}^n$, where $n \in \mathbb{Z}_{\geq 0}$ is a suitable integer. Let $M$ be the dual lattice of $N$. We then consider a rational polyhedral cone $\sigma \subseteq N \otimes_{\mathbb{R}} \mathbb{R} \cong \mathbb{R}^n$. To this cone, we associate the semigroup $S_\sigma = \sigma^\vee \cap M$. The corresponding affine toric variety $U_\sigma$ is given by\cite{CLS11}:
\begin{equation}
U_\sigma = \mathrm{Spec} \left(\mathbb{C} \left[ S_\sigma \right] \right) = \mathrm{Spec} \left(\mathbb{C} \left[ \sigma^\vee \cap M \right]\right) \, .
\end{equation}
As an example, consider
\begin{equation}
\sigma = \mathrm{Span}_{\mathbb{Z}_{\geq 0}} \left( \left[ \begin{array}{c} 1 \\ 0 \end{array} \right], \left[ \begin{array}{c} 0 \\ 1 \end{array} \right] \right) \, .
\end{equation}
We create $U_\sigma$ in \texttt{OSCAR}:
\begin{minted}{Julia}
o = positive_hull([1 0; 0 1])
U = affine_normal_toric_variety(o)
\end{minted}
Many properties of $U_\sigma$ are encoded in $\sigma$. For instance, $U_\sigma$ is smooth if and only if $\sigma$ can be generated by a subset of a basis of the lattice $N$. An interactive check in \texttt{OSCAR} can determine whether $U_\sigma$ is smooth:
\begin{minted}{Julia}
julia> hilbert_basis(o)
2-element SubObjectIterator
{PointVector{ZZRingElem}}:
 [1, 0]
 [0, 1]

julia> is_smooth(U)
true
\end{minted}
Similarly, the dimension of $U_\sigma$ matches that of $\sigma$:
\begin{minted}{Julia}
julia> dim(o) == dim(U)
true
\end{minted}
Below is an instance of a non-smooth affine toric variety that can be created using \texttt{OSCAR}:
\begin{minted}{Julia}
o2=positive_hull([-1 1; 0 1; 1 1])
U2=affine_normal_toric_variety(o2)
\end{minted}
We verify interactively that $U_2$ is not smooth:
\begin{minted}{Julia}
julia> hilbert_basis(o2)
3-element SubObjectIterator
{PointVector{ZZRingElem}}:
 [-1, 1]
 [0, 1]
 [1, 1]

julia> is_smooth(U2)
false
\end{minted}
Notice the appearance of the generator $[0, 1]$ in the Hilbert basis of $\sigma_2$. Its appearance signifies that $U_2$ is not smooth. Alternatively, we can inspect $U_2$ as subvariety of the affine space. To this end, we compute the toric ideal:
\begin{minted}{Julia}
julia> toric_ideal(U2)
ideal(-x1*x2 + x3^2)
\end{minted}
This means, that in the affine space $\mathbb{A}^3$ with coordinates $(x_1, x_2, x_3)$, it holds $U_{2} \cong V( - x_1 x_2 + x_3^2 )$. Consequently, $U_2$ is singular.

Much more can be said about the interplay between polyhedral and toric geometry. For example, there exists a connection between normal toric varieties and rational polyhedral fans. Indeed, in \texttt{OSCAR}, one can create a general normal toric variety based on a rational polyhedral fan. Moreover, \texttt{OSCAR} offers specialized constructors for several well-known toric varieties, including \texttt{projective$\_$space}, \texttt{hirzebruch$\_$surface}, \texttt{del$\_$pezzo$\_$surface}, and \texttt{cyclic$\_$quotient$\_$singularity}. For further information, interested readers may wish to consult \cite{CLS11}.

\Ueberschrift{Notable capabilities}{sec3}

\Ueberschriftu{Vanishing sets of line bundle cohomology}
Support for torus invariant divisors, divisor classes and line bundles is available in \texttt{OSCAR}. The \texttt{cohomCalg} algorithm \cite{BJRR10, **key**2} is employed to infer dimensions of line bundle cohomologies on any smooth and complete, as well as any simplicial and projective toric variety $X_\Sigma$. The set $V^i(X_\Sigma)$ of all line bundles on $X_\Sigma$ with vanishing $i$-th sheaf cohomology can be derived \cite{Bie18}:
\begin{align}
V^i (X_\Sigma) = \mathrm{Pic} (X_\Sigma) - \bigcup_{m = 1}^{l}{L^i_{(m)}} \, ,
\end{align}
where $L^i_{(m)}$ is the set of lattice points in a certain polyhedron $P^i_{(m)}$. For $\mathbb{P}^1 \times \mathbb{P}^1$, it holds $\mathrm{Pic}( \mathbb{P}^1 \times \mathbb{P}^1) = \mathbb{Z}^2$ and that the vanishing sets can be represented as follows:
\begin{equation}
\begin{tikzpicture}[baseline=(current  bounding  box.center)]

    \def\s{0.8};
    
    \draw { (0,0)--(3*\s,0)--(3*\s,3*\s)--(0,3*\s)--(0,0) } [opacity = 0.5, fill=green];
    \draw { (-2*\s,0)--(-5*\s,0)--(-5*\s,3*\s)--(-2*\s,3*\s)--(-2*\s,0) } [opacity = 0.5, fill=yellow];
    \draw { (0,-2*\s)--(3*\s,-2*\s)--(3*\s,-5*\s)--(0,-5*\s)--(0,-2*\s) } [opacity = 0.5, fill=red];
    \draw { (-2*\s,-2*\s)--(-2*\s,-5*\s)--(-5*\s,-5*\s)--(-5*\s,-2*\s)--(-2*\s,-2*\s) } [opacity = 0.5, fill=blue];

    \node at (1.5*\s,1.5*\s) {$P^0$};
    \node at (-3.5*\s,1.5*\s) {$P^1_{(1)}$};
    \node at (1.5*\s,-3.5*\s) {$P^1_{(2)}$};
    \node at (-3.5*\s,-3.5*\s) {$P^2$};

    \draw (2*\s,-0.1) -- (2*\s,0.1);
    \draw (1*\s,-0.1) -- (1*\s,0.1);
    \draw (-1*\s,-0.1) -- (-1*\s,0.1);
    \draw (-2*\s,-0.1) -- (-2*\s,0.1);
    \draw (-3*\s,-0.1) -- (-3*\s,0.1);
    \draw (-4*\s,-0.1) -- (-4*\s,0.1);
    \draw (-0.1,2*\s) -- (0.1,2*\s);
    \draw (-0.1,1*\s) -- (0.1,1*\s);
    \draw (-0.1,-1*\s) -- (0.1,-1*\s);
    \draw (-0.1,-2*\s) -- (0.1,-2*\s);
    \draw (-0.1,-3*\s) -- (0.1,-3*\s);
    \draw (-0.1,-4*\s) -- (0.1,-4*\s);

    \node at (1*\s,-0.35) {$1$};
    \node at (-0.25, 1*\s) {$1$};

    \draw [thick,-latex] (0,-5*\s)--(0,3*\s) node (yaxis) [above] {};
    \draw [thick,-latex] (-5*\s,0)--(3*\s,0) node (xaxis) [right] {};
    
\end{tikzpicture}
\label{equ:Picture1}
\end{equation}
Specifically,
\begin{align}
\begin{split}
&V^0(\mathbb{P}^1 \times \mathbb{P}^1) = \mathbb{Z}^2 - (P^0 \cap \mathbb{Z}^2) \, ,\\
&V^1(\mathbb{P}^1 \times \mathbb{P}^1) = \mathbb{Z}^2 - (P^1_{(1)} \cup P^1_{(2)}) \cap  \mathbb{Z}^2\, , \\
&V^2(\mathbb{P}^1 \times \mathbb{P}^1) = \mathbb{Z}^2 - (P^2 \cap \mathbb{Z}^2) \, , \\
&P^0 = \left[ \begin{array}{c} 0 \\ 0 \end{array} \right] + \mathrm{Span}_{\mathbb{Z}_{\geq 0}} \left( \left[ \begin{array}{c} 1 \\ 0 \end{array} \right], \left[ \begin{array}{c} 0 \\ 1 \end{array} \right] \right) \, , \\
&P^1_{(1)} = - \left[ \begin{array}{c} 2 \\ 0 \end{array} \right] + \mathrm{Span}_{\mathbb{Z}_{\geq 0}} \left( \left[ \begin{array}{c} -1 \\ 0 \end{array} \right], \left[ \begin{array}{c} 0 \\ 1 \end{array} \right] \right) \, ,\\
&P^1_{(2)} = - \left[ \begin{array}{c} 0 \\ 2 \end{array} \right] + \mathrm{Span}_{\mathbb{Z}_{\geq 0}} \left( \left[ \begin{array}{c} 1 \\ 0 \end{array} \right], \left[ \begin{array}{c} 0 \\ -1 \end{array} \right] \right) \, , \\
&P^2 = - \left[ \begin{array}{c} 2 \\ 2 \end{array} \right] - \mathrm{Span}_{\mathbb{Z}_{\geq 0}} \left( \left[ \begin{array}{c} 1 \\ 0 \end{array} \right], \left[ \begin{array}{c} 0 \\ 1 \end{array} \right] \right) \, .
\end{split}
\label{equ:VanishingSetsExplicitly}
\end{align}
With the following lines, we can replicate these results in \texttt{OSCAR}:

\begin{minted}{Julia}
P1 = projective_space(
            NormalToricVariety, 1)
v0, v1, v2 = vanishing_sets(P1*P1)
ph0 = polyhedra(v0)[1]
ph11, ph12 = polyhedra(v1)
ph2 = polyhedra(v2)[1]
\end{minted}
With \texttt{OSCAR}, we can investigate the polyhedra interactively. For example, we can find inequalities for $P^0$ and $P^2$ as follows:
\begin{minted}{Julia}
julia> print_constraints(ph0)
-x₁ ≦ 0
-x₂ ≦ 0

julia> print_constraints(ph2)
x2 ≦ -2
x1 ≦ -2
\end{minted}
Indeed, from \cref{equ:VanishingSetsExplicitly} we see that $L^0$, $L^2$ can be expressed as follows:
\begin{align}
L^0 &= \{ \left. (x_1,x_2) \in \mathbb{Z}^2 \right|  x_1, x_2 \geq 0 \} \, , \\
L^2 &= \{ \left. (x_1,x_2) \in \mathbb{Z}^2 \right|  x_1, x_2 \leq -2 \} \, .
\end{align}
It is not too hard to repreat this exercise for $L^1_{(1)}$, $L^1_{(2)}$.

As a more interesting example, consider the del Pezzo surface $dP_1$ with $\mathbb{Z}^2$-graded Cox ring:
\begin{align}
\begin{tabular}{c|ccc|c}
\toprule
& $x_1$ & $x_2$ & $x_3$ & $e_1$ \\
\midrule
$H$     & $1$ & $1$ & $1$ & \\
$-E_1$ & $1$ & $1$ &    & $-1$ \\
\bottomrule
\end{tabular}
\label{equ:GradingCoxRingDP1}
\end{align}
For this grading, we visualize the vanishing sets:
\begin{equation}
 \begin{tikzpicture}[baseline=(current  bounding  box.center)]

    \def\s{0.8};
    
    \draw { (0,0)--(3*\s,3*\s)--(3*\s,-5*\s)--(0,-5*\s)--(0,0) } [opacity = 0.5, fill=green];
    \draw { (-2*\s,-2*\s)--(3*\s,-2*\s)--(3*\s,-5*\s)--(-5*\s,-5*\s)--(-2*\s,-2*\s) } [opacity = 0.5, fill=red];
    \draw { (-1*\s,1*\s)--(1*\s,3*\s)--(-5*\s,3*\s)--(-5*\s,1*\s)--(-1*\s,1*\s) } [opacity = 0.5, fill=yellow];
    \draw { (-3*\s,-1*\s)--(-3*\s,3*\s)--(-5*\s,3*\s)--(-5*\s,-3*\s)--(-3*\s,-1*\s) } [opacity = 0.5, fill=blue];

    \node at (1.5*\s,-1.0*\s) {$P^0$};
    \node at (-1.5*\s,2.0*\s) {$P^1_{(1)}$};
    \node at (-1.5*\s,-3.5*\s) {$P^1_{(2)}$};
    \node at (-4*\s,-1.0*\s) {$P^2$};
    
    \draw (2*\s,-0.1) -- (2*\s,0.1);
    \draw (1*\s,-0.1) -- (1*\s,0.1);
    \draw (-1*\s,-0.1) -- (-1*\s,0.1);
    \draw (-2*\s,-0.1) -- (-2*\s,0.1);
    \draw (-3*\s,-0.1) -- (-3*\s,0.1);
    \draw (-4*\s,-0.1) -- (-4*\s,0.1);
    \draw (-0.1,2*\s) -- (0.1,2*\s);
    \draw (-0.1,1*\s) -- (0.1,1*\s);
    \draw (-0.1,-1*\s) -- (0.1,-1*\s);
    \draw (-0.1,-2*\s) -- (0.1,-2*\s);
    \draw (-0.1,-3*\s) -- (0.1,-3*\s);
    \draw (-0.1,-4*\s) -- (0.1,-4*\s);
        
    \node at (1*\s,-0.35) {$1$};
    \node at (-0.25, 1*\s) {$1$};
    
    \draw [thick,->] (0,-5*\s)--(0,3*\s) node (yaxis) [above] {};
    \draw [thick,->] (-5*\s,0)--(3*\s,0) node (xaxis) [right] {};

\end{tikzpicture}
\label{equ:Picture2}
\end{equation}
The interested reader might find it entertaining to ``see'' Serre duality in \cref{equ:Picture1} and \cref{equ:Picture2}.

We emphasize that the vanishing sets can be determined algorithmically for any smooth and complete, as well as any simplicial and projective toric variety $X_\Sigma$. However, our ability to visualize the vanishing sets reduces drastically once the polyhedra are of dimension 4 or higher. This happens for instance for $\mathbb{P}^1 \times \mathbb{P}^1 \times \mathrm{dP}_1$. Still, the vanishing sets can be derived in \texttt{OSCAR}:
\begin{minted}{Julia}
P1 = projective_space(
             NormalToricVariety,1)
dP1 = del_pezzo_surface(1)
v0, v1, v2, v3, v4
       = vanishing_sets(P1*P1*dP1)
\end{minted}

\Ueberschriftu{Intersection theory}
Loosely speaking, intersection theory provides an answer to the question ``At how many points do two algebraic cycles intersect?''`. A caveat arises whenever the algebraic cycles in question are ``similar/the same''. This leads to the notion of \emph{rational equivalence} and the observation, that sometimes the number of intersection points can be negative. To demonstrate this somewhat exotic idea in a concrete setting, we focus on the del Pezzo surface $dP_1$ with $\mathbb{Z}^2$-graded Cox ring as in \cref{equ:GradingCoxRingDP1}. Next, consider the following algebraic cycles:
\begin{align}
H = V(x_1) + V(e_1) \, , \qquad E_1 = V(e_1) \, .
\end{align}
Strictly speaking, we want to consider the rational equivalence classes of these algebraic cycles. For ease of notation, we do not introduce new symbols. The intersection numbers among $H$ and $E_1$ are as follows:
\begin{align}
H^2 = 1 \, , \qquad H \cdot E_1 = 0 \, , \qquad E_1 \cdot E_1 = - 1 \, .
\end{align}
The following code computes this result in \texttt{OSCAR}:
\begin{minted}{Julia}
julia> dP1 = del_pezzo_surface(1);

julia> intersection_form(dP1)
Dict{MPolyRingElem, QQFieldElem}
with 10 entries:
  x1*x3 => 1
  e1^2  => -1
  x2*x3 => 1
  x3^2  => 1
  x1*x2 => 0
  x3*e1 => 0
  x2^2  => 0
  x1*e1 => 1
  x2*e1 => 1
  x1^2  => 0
\end{minted}
Certainly, we can create the rational equivalence classes of $H$ and $E_1$ in \texttt{OSCAR}:
\begin{minted}{Julia}
x1, x2, x3, e1
            = gens(chow_ring(dP1))
E1 = 
rational_equivalence_class(dP1,e1)
H = E1 +
rational_equivalence_class(dP1,x1)
\end{minted}
With this, we can explicitly and interactively verify in \texttt{OCSAR} how these algebraic cycles intersect:
\begin{minted}{Julia}
julia> H*H
Rational equivalence class on a
normal toric variety represented
by V(x2,x3)

julia> H*E1
Trivial rational equivalence class
on a normal toric variety

julia> E1*E1
Rational equivalence class on a
normal toric variety represented
by -1V(x2,x3)
\end{minted}
In the last computation, notice the appearance of $-1$. To understand its meaning, we must understand how the intersection points of $E_1$ with itself are computed. The theory tells us that we should use different, yet rationally equivalent, algebraic cycles which intersect ``nicely''. The technical term for this is to move the algebraic cycles in \emph{general position} \cite{EH16}.

In toric varieties, rational equivalences are captured by the ideal of linear relations:
\begin{minted}{Julia}
julia> ideal_of_linear_\
                    relations(dP1)
ideal(x1 - x3 + e1, x2 - x3 + e1)
\end{minted}
Let $\sim$ denote rational equivalence. Then it holds:
\begin{align}
V(x_1) - V(x_3) + V(e_1) &\sim 0 \, , \\
V(x_2) - V(x_3) + V(e_1) &\sim 0 \, .
\end{align}
Hence $V(e_1) \sim V(x_3) - V(x_1)$ and it follows that
\begin{align}
E_1 \cdot E_1 &\sim V(e_1) \cdot \left[ V(x_3) - V(x_1) \right] \\
                      &= V( e_1, x_3) - V( e_1, x_1) \, .
\end{align}
Next, let us look at the Stanley-Reisner ideal of $\mathrm{dP}_1$:
\begin{minted}{Julia}
julia> stanley_reisner_ideal(dP1)
ideal(x1*x2, x3*e1)
\end{minted}
From this ideal we learn that
\begin{align}
\left\{ \left. p \in \mathrm{dP}_1 \right| x_3 = e_1 = 0 \right\} = \emptyset \, .
\end{align}
Consequently, we find
\begin{align}
E_1 \cdot E_1 \sim - V( e_1, x_1) \, .
\end{align}
It is not too hard to verify that $V(e_1, x_1) \sim V(x_2, x_3)$. This finally aligns our investigations with the result computed by \texttt{OSCAR}. We do hope that this example illustrates the origin of negative intersection numbers.

The collection of rational equivalence classes of algebraic cycles enjoys a ring structure where the multiplication corresponds to the intersection of the algebraic cycles. This ring is known as the \emph{Chow ring} and can be computed for any complete, simplicial toric variety \cite{CLS11}. For example, the Chow ring of $dP_1$ can be computed in \texttt{OSCAR} as follows:
\begin{minted}{Julia}
julia> chow_ring(dP1)
Quotient of Multivariate
Polynomial Ring in x1, x2, x3, e1
over Rational Field by ideal
(x1-x3+e1, x2-x3+e1, x1*x2, x3*e1)
\end{minted}
It has been noted more recently that the completeness assumption can be dropped \cite{Peg14}.\footnote{See also \cite{FY04} for the significance of this observation for the interplay between matroids and toric varieties.} Indeed, \texttt{OSCAR} is capable of computing the Chow ring for simplicial toric varieties that are not complete. As an example, we create a non-complete, yet simplicial toric variety \texttt{v} from its rays \texttt{r} and (maximal) cones \texttt{c}:
\begin{minted}{Julia}
r =[[1,0],[0,1],[-1,-1]]
c = [[1],[2],[3]]
v = normal_toric_variety(r, c)
\end{minted}
We verify that $v$ is not complete but simplicial:
\begin{minted}{Julia}
julia> is_complete(v)
false

julia> is_simplicial(v)
true
\end{minted}
We can also compute the Chow ring interactively:
\begin{minted}{Julia}
julia> chow_ring(v)
Quotient of Multivariate
Polynomial Ring in x1, x2, x3
over Rational Field by ideal
(x1-x3, x2-x3, x1*x2, x1*x3, x2*x3)
\end{minted}

\Ueberschriftu{Triangulations}
As explained in \cite{KS98, KS00}, reflexive polytopes $\Delta^\circ$ (and their polar duals) can be used to classify Calabi-Yau hypersurfaces in toric spaces. The ambient toric spaces can be found from fine regular star triangulations (FRST) of $\Delta^\circ$ (see e.g. \cite{DRS10} for background). For an example, consider the square with vertices at $(\pm 1, \pm 1)$. The informed reader will notice immediately that this configuration has a unique FRST corresponding to $\mathbb{P}^1 \times \mathbb{P}^1$.
\begin{minted}{Julia}
P = convex_hull(
         [1 1; -1 1; 1 -1; -1 -1])
X = NormalToricVarieties\
         FromStarTriangulations(P)
 \end{minted}
Certainly, we can verify that $X$ consists only of a single variety. Furthermore, we compute the Stanley-Reisner and the irrelevant ideal of this toric variety, to provide evidence that this variety is indeed just 
$\mathbb{P}^1 \times \mathbb{P}^1$:
\begin{minted}{Julia}
julia> length(X)
1

julia> irrelevant_ideal(X[1])
ideal(x3*x4, x2*x4, x1*x3, x1*x2)

julia> stanley_reisner_ideal(X[1])
ideal(x1*x4, x2*x3)
\end{minted}
A much more involved example is included in the tutorial (\href{https://www.oscar-system.org/tutorials/}{https://www.oscar-system.org/tutorials/}). This example is computationally demanding and its code was optimized for performance. We propose to use this example for benchmarking purposes. Note also that this code was used in a recent string theory application \cite{BCDO22}.

\Ueberschrift{Outlook}{sec5}

\texttt{OSCAR} is a relatively new software and still under heavy development. This is an opportunity for young developers -- we truly appreciate contributions. While this means that things are changing within \texttt{OSCAR}, the interface for toric varieties is already rather mature. In recent times, this interface has remained stable. For users (of the toric functionality) this is great news, as you need not be afraid of changes to the interface that might break your workflow. We strongly encourage users to try out and enjoy the existing toric functionality.

There are plans to significantly extend the toric functionality for coherent sheaves. In the realm of smooth and complete toric varieties, coherent sheaves are equivalent to certain classes of finitely presented graded modules \cite{CLS11}. This equivalence can be utilized to compute sheaf cohomologies of coherent sheaves. In fact, a relevant algorithm for this purpose was proposed in \cite{Bie18}.\footnote{This algorithm is available at \href{https://github.com/homalg-project/ToricVarieties_project}{https://github.com/homalg-project/ToricVarieties$\_$project}.} There are plans to incorporate this functionality into \texttt{OSCAR} in the future. It would also be advantageous to explore specialized algorithms, e.g. based on \cite{Pay08, ABKW18, AP19, AW21}, for cohomologies of vector bundles.

In view of applications in the field of F-theory -- a specialized domain of string theory -- initial discussions have taken place to assess the possibility of incorporating \texttt{FTheoryTools} \cite{BT22} into \texttt{OSCAR}. \texttt{FTheoryTools} is primarily focused on computing resolution for singular elliptic fibrations. Such computations pose a significant arithmetic challenge in F-theory. The goal is to make this task as convenient as possible for researchers in this area. There are overlaps with some of the schemes technology that is currently being actively developed in \texttt{OSCAR}. For instance, \texttt{OSCAR} has basic support for toric schemes in experimental stages.

A more specialized task in F-theory involves constructing solutions that replicate the particle physics observed in modern accelerator experiments. Recently, numerous promising solutions known as the \emph{Quadrillion F-theory Standard Models} (F-theory QSMs) were identified in \cite{CHLLT19}. These solutions are based on the geometry of toric K3 surfaces via \cite{KS98}. Consequently, toric technology is critical to constructing and exploring these solutions in the future. In fact, many F-theory constructions are based on toric geometry (see \cite{Wei18} and references therein). It would be interesting to provide user-friendly and convenient tools in \texttt{OSCAR} for toric F-theory constructions.

Cosmological investigations within string theory led to the development of the software \emph{CYTools} \cite{DRM22}. This software focuses on high-performance triangulations of the 4-dimensional reflexive polytopes in \cite{KS00}. Such triangulation tasks matter in many explicit realizations of Calabi-Yau manifolds. For these reasons, \cite{DRM22} is a very interesting software package. We expect that its capabilities can be boosted by using \texttt{mptopcom} \cite{JJK18} or the latest version of \texttt{TOPCOM} \cite{Ram02}. This task is reserved for future work.

\subsection*{Acknowledgement}
M.~B. and L.~K. express their gratitude and appreciation for the support provided by the OSCAR team, led by Claus Fieker, Max Horn, Michael Joswig, and Wolfram Decker. M.~B. acknowledges financial support from the Forschunginitiative des Landes Rheinland-Pfalz through the project \emph{SymbTools – Symbolic Tools in Mathematics and their Application}. L. K. is thankful for the funding received from \emph{MaRDI -- Mathematical research initiative} of the German Research Foundation (DFG). This work was supported by the SFB-TRR 195 \emph{Symbolic Tools in Mathematics and their Application} of the German Research Foundation (DFG).

\bibliographystyle{abbrv}
\bibliography{references}{}

\begin{thebibliography}{10}

\bibitem{**key**2}
{cohomCalg package -- High-performance line bundle cohomology computation based
  on \cite{BJRR10}}.
\newblock {https://github.com/BenjaminJurke/cohomCalg}, 2010.

\bibitem{ABKW18}
K.~Altmann, J.~Buczyński, L.~Kastner, and A.-L. Winz.
\newblock {Immaculate line bundles on toric varieties}, 2018.

\bibitem{AP19}
K.~Altmann and D.~Ploog.
\newblock {Displaying the cohomology of toric line bundles}, 2019.

\bibitem{AW21}
K.~Altmann and F.~Witt.
\newblock {The structure of exceptional sequences on toric varieties of Picard
  rank two}, 2021.

\bibitem{BEKS17}
J.~Bezanson, A.~Edelman, S.~Karpinski, and V.~B. Shah.
\newblock {Julia: A fresh approach to numerical computing}.
\newblock {\em SIAM {R}eview}, 59(1):65--98, 2017.

\bibitem{Bie18}
M.~Bies.
\newblock {\em {Cohomologies of coherent sheaves and massless spectra in
  F-theory}}.
\newblock PhD thesis, Heidelberg U., 2 2018.

\bibitem{BCDO22}
M.~Bies, M.~Cveti\v{c}, R.~Donagi, and M.~Ong.
\newblock {Brill-Noether-general limit root bundles: absence of vector-like
  exotics in F-theory Standard Models}.
\newblock {\em JHEP}, 11:4, 2022.

\bibitem{BT22}
M.~Bies and A.~P. Turner.
\newblock {FTheoryTools -- Julia tools for F-Theory compactifications}.
\newblock \url{https://github.com/Julia-meets-String-Theory/FTheoryTools.jl},
  2022-2023.

\bibitem{BJRR10}
R.~Blumenhagen, B.~Jurke, T.~Rahn, and H.~Roschy.
\newblock {Cohomology of Line Bundles: A Computational Algorithm}.
\newblock {\em J. Math. Phys.}, 51:103525, 2010.

\bibitem{CLS11}
D.~A. Cox, J.~B. Little, and H.~K. Schenck.
\newblock {\em {Toric Varieties}}.
\newblock Graduate studies in mathematics. American Mathematical Soc., 2011.

\bibitem{CHLLT19}
M.~Cveti{\v c}, J.~Halverson, L.~Lin, M.~Liu, and J.~Tian.
\newblock {Quadrillion $F$-Theory Compactifications with the Exact Chiral
  Spectrum of the Standard Model}.
\newblock {\em Phys. Rev. Lett.}, 123(10):101601, 2019.

\bibitem{DRS10}
J.~A. {De Loera}, J.~Rambau, and F.~Santos.
\newblock {\em {Triangulations}}.
\newblock Graduate studies in mathematics. Springer, Heidelberg Dordrecht
  London New York, 2010.

\bibitem{DRM22}
M.~Demirtas, A.~Rios-Tascon, and L.~McAllister.
\newblock {CYTools: A Software Package for Analyzing Calabi-Yau Manifolds},
  2022.

\bibitem{**key*}
{Eder, C. and Decker, W. and Fieker, C. and Horn, M. and Joswig, M.}, editor.
\newblock {\em {The OSCAR book}}.
\newblock {2024}.

\bibitem{EH16}
D.~Eisenbud and J.~Harris.
\newblock {\em {3264 and All That: A Second Course in Algebraic Geometry},}.
\newblock Cambridge University Press, 2016.

\bibitem{FY04}
E.~M. Feichtner and S.~Yuzvinsky.
\newblock Chow rings of toric varieties defined by atomic lattices.
\newblock {\em Inventiones Mathematicae}, 155(3):515--536, mar 2004.

\bibitem{GS23}
D.~R. Grayson and M.~E. Stillman.
\newblock {Macaulay2, a software system for research in algebraic geometry}.
\newblock Available at \url{http://www.math.uiuc.edu/Macaulay2/}, 2023.

\bibitem{JJK18}
C.~Jordan, M.~Joswig, and L.~Kastner.
\newblock {Parallel Enumeration of Triangulations}.
\newblock {\em The Electronic Journal of Combinatorics}, 2018.
\newblock https://polymake.org/doku.php/mptopcom.

\bibitem{KS98}
M.~Kreuzer and H.~Skarke.
\newblock {Classification of reflexive polyhedra in three-dimensions}.
\newblock {\em Adv. Theor. Math. Phys.}, 2:853--871, 1998.

\bibitem{KS00}
M.~Kreuzer and H.~Skarke.
\newblock {Complete classification of reflexive polyhedra in four-dimensions}.
\newblock {\em Adv. Theor. Math. Phys.}, 4:1209--1230, 2000.

\bibitem{**key**3}
{OSCAR -- Open Source Computer Algebra Research system, Version 0.12.0-DEV}.
\newblock \url{https://www.oscar-system.org}, 2023.

\bibitem{Pay08}
S.~Payne.
\newblock {Moduli of toric vector bundles}.
\newblock {\em Compositio Mathematica}, 144(5):1199--1213, Sept 2008.

\bibitem{Peg14}
C.~Pegel.
\newblock Chow rings of toric varieties.
\newblock Master's thesis, University of Bremen, Faculty of Mathematics, Sept
  2014.
\newblock Refereed by Prof. Dr. Eva Maria Feichtner and Dr. Emanuele Delucchi.

\bibitem{Ram02}
J.~Rambau.
\newblock {TOPCOM: Triangulations of Point Configurations and Oriented
  Matroids}.
\newblock {\em Proceedings of the International Congress of Mathematical
  Software}, 2002.

\bibitem{Dev23}
{{S}ageMath, the {S}age {M}athematics {S}oftware {S}ystem ({V}ersion 9.8)}.
\newblock \url{https://www.sagemath.org}, 2023.

\bibitem{Wei18}
T.~Weigand.
\newblock {F-theory}.
\newblock {\em PoS}, TASI2017:016, 2018.

\end{thebibliography}

\end{multicols}

\end{otherlanguage}

\end{document}